\documentclass[a4paper,12pt]{amsart} 
\usepackage{amssymb,euscript,eufrak,verbatim}
\usepackage[utf8]{inputenc}      
\usepackage[T1]{fontenc}
\usepackage[polish,english]{babel}

\usepackage{tikz}
\usetikzlibrary{snakes}
\usetikzlibrary{arrows}

\def\z{\phi}
\def\e{\varepsilon}

\def\U{\mathcal{U}}

\def\N{\mathbb{N}}
 
\def\S1{\mathbb{S}^1}

\def\R{\mathbb{R}}

\def\B{\mathcal{B}}

\newcommand{\vt}{{\vec\theta}}
\newcommand{\zt}{\phi^{\vt}}

\def\G{\mathcal{G}}
\def\i{{i,j}}

\newtheorem{theorem}{Theorem}

\newtheorem{proposition}[theorem]{Proposition}
\newtheorem{corollary}[theorem]{Corollary}
\newtheorem*{remark}{Remark}
\DeclareMathSymbol{\varnothing}{\mathord}{AMSb}{"3F} 
\begin{document}
\renewenvironment{proof}{\noindent {\bf Proof.}}{ \hfill\qed\\ }
\newenvironment{proofof}[1]{\noindent {\bf Proof of #1.}}{ \hfill\qed\\ }

\title{Infinite ergodic index of the Ehrenfest wind-tree model}  

\author {Alba M\'alaga Sabogal}
\address{CHArt/THIM (EA 4004)
Université Paris 8-Vincennes-Saint-Denis.}
\email{alba.malaga@polytechnique.edu}

\def\curraddrname{{\itshape Address}}
\author{Serge Troubetzkoy}
\address{Aix Marseille Univ, CNRS, Centrale Marseille, I2M,  Marseille, France}

  \curraddr{I2M, Luminy\\ Case 907\\ F-13288 Marseille CEDEX 9\\ France}

 \email{serge.troubetzkoy@univ-amu.fr}
 
 \thanks{We thank Jack Milnor for suggesting a nice presentation of our topology.
 AMS acknowledges that this work was started during a post-doc funded by the A*MIDEX project (ANR-11-IDEX-0001-02), funded itself by the “Investissements d’avenir” program of the French Government, managed by the French National Research Agency (ANR)”. She continued working on this project  during the ATER position she held at Mathematics Laboratory in Orsay in 2015-2016. ST gratefully acknowledges the support of project  APEX "Systèmes dynamiques: Probabilités et Approximation Diophantienne PAD" funded by the Région PACA}
\begin{abstract}   
The set of all possible configurations of the Ehrenfest wind-tree model endowed with the Hausdorff topology  is a compact metric space. For a typical configuration we show that the wind-tree dynamics has infinite ergodic index in almost every direction.
In particular some ergodic theorems can be applied to show that if we start with a large number of initially parallel
particles their directions decorrelate as the dynamics evolve answering the question posed by the Ehrenfests.
\end{abstract} 
\maketitle

\section{Introduction}
In 1912 Paul and Tatiana Ehrenfest wrote a seminal article on the foundations of Statistical Mechanics in which  the wind-tree model was introduced
 in order to interpret the work of Boltzmann and Maxwell on gas dynamics \cite{EhEh}. 
In the wind-tree model a point particle moves without friction on the plane with infinitely many rigid obstacles removed, and collides elastically with the obstacles.	
 The Ehrenfests' paper dates from times when the notions of probability theory where not yet rigorously defined. Thus they could not describe the  distribution of the obstacles in a probabilistic way, they used the word ``irregular'' to describe it. However, they made precise what they did expected from the placement of the obstacles: obstacles are identical squares, all parallel to each other,
 the placement is  irregular,
 every portion of the plane contains about the same number of obstacles,
and the distances between the obstacles are large in comparison to the obstacle's size.

If we fix the direction of the particle, the billiard flow will take only four directions.
The Ehrenfests asked the following question: start $K$ particles in a given direction, will the number of
particles in each of the four directions asymptotically equalize to about $K/4$?
To answer this question we study the ergodic properties of the wind-tree model. 
Interestingly the birth of ergodic theory can be traced back to the Ehrenfests' article in which the word ergodic was used for
the first time  with a  close mathematical meaning to the current one  \cite{GaBoGe}. 
We  consider the set of all possible configurations and introduce a canonical topology  which makes it
 a compact metric space.  We show that for (Baire) generic  configurations, for almost every direction the billiard flow has infinite ergodic index, i.e., all its powers are ergodic. In a finite measure space this would be equivalent to saying that the
 flow is weakly mixing.   The asymptotic equalization of the directions of $K$ particles  in several senses
 then follows from various ergodic theorems (note that we are in the framework of infinite ergodic theory here, so the Birkhoff ergodic theorem is not directly applicable).

	In two previous articles we have considered a subset of configurations which are small perturbations of lattice configurations, and we showed  that the generic wind-tree is minimal and ergodic in  almost every direction \cite{MSTr1,MSTr2}.  The topology considered
in these two articles is equivalent to the one considered in this article.  Furthermore the proofs of these two results
hold mutatis  mutandis  in the more general setting 	which we consider here.
	
	There have been a number of  results on the wind-tree model \cite{DeCoVB, Ga, HaCo, HaCo1, Tr1, VBHa, WoLa}, and  on the wind-tree model
	 with periodical distribution of obstacles of squares, rectangles and more recently other polygonal shapes \cite{AvHu,BaKhMaPl,BiRo, De,DeHuLe,DeZo,FrHu,FrUl,HaWe,HuLeTr}.

\section{Definitions and main results}

 For sake of simplicity, a square whose sides are parallel to lines $y=\pm x$ will be referred to as \emph{rhombus} in the rest of the article. The $\mathcal L_1$ distance in $\mathbb R^2$ will be denoted by $d$. Note that balls with respect to this distance are rhombii.

Fix $s > 0$. A \emph{configuration} is an at most countable collection of rhombii with diameter $s$, whose interiors are pairwise disjoint. 
Since $s$ is fixed it is enough to note the centers of the rhombii, thus 
a configuration $g$ is  an at most countable subset of $\mathbb R^2$ such that if $z_1,z_2 \in g$ then $d(z_1,z_2) \ge s$.

To define a topology on the set of configurations consider polar coordinates $(r,\theta)$ on the plane. Each point $(r,\theta)$ in the plane is the stereographic projection of a point in the sphere with spherical coordinates $(2\arctan(1/r),\theta)$.   
Apply  the inverse of the stereographic projection to a configuration $g$ to obtain a subset of the sphere.  
Let $\hat{g}$ denote the union of this  set with the north pole of the sphere denoted by  $\{\infty\}$, it is a closed subset of the sphere.
The topology we define on the set of configurations is then induced by the Hausdorff distance $d_H$ 
given by  
$$d_H(g_1,g_2) = \max(\sup\limits_{z_1\in \hat{g}_1} \inf\limits_{z_2\in \hat{g}_2} \rho(z_1,z_2) , \sup\limits_{z_2\in \hat{g}_2} \inf\limits_{z_1\in \hat{g}_1} \rho(z_1,z_2)).$$
Here $\rho$ denotes the geodesic distance on the sphere, i.e., the length of the shortest path from one point to another along the great circle passing through them.  
Let $\textit{Conf}$ be the set of all configurations.

\begin{proposition}\label{p1}
$(\textit{Conf},d_H)$ is compact metric space, thus a Baire space.
\end{proposition}

The proposition is proven in the appendix.
Let $\U_\varepsilon(g)$ be the set of all configurations that are at most $\varepsilon$-close to $g$
$$\U_\varepsilon(g) := \{g' | d_H(g',g)<\varepsilon\}.$$

\begin{proposition}
There is a dense $G_{\delta}$ subset $G$ of $(\textit{Conf},d_H)$ such that for each $g \in G$ 
\begin{enumerate}
\item $g$ is an infinite configuration,
\item every pair of points $z_1,z_2 \in g$ satisfy $d(z_1,z_2) > s$.
\end{enumerate}
\end{proposition}

\begin{remark}
Point (2) means that the obstacles centered at $z_1$ and $z_2$ do not intersect.
\end{remark}

\begin{proof}
There are infinite configurations arbitrarily close to any finite configuration,  thus we can choose a countable dense set  $\{g_n : n \in \mathbb{N}\}$ of
infinite configurations. 
Let $\varepsilon(g_n)> 0$  be the infinimum of $\varepsilon$ such  that there are at least $n$ distinct points $z \in g_n$ satisfying $\rho(z,\infty) > \varepsilon$. Denote this necessarily finite set of points by $B(g_n)$.
Clearly
$$G := \bigcap_{m=1}^\infty \bigcup_{n \ge m} U_{\varepsilon(g_n)} (g_n)$$
is a dense $G_\delta$ set.
If $g \in G$, then $g$ is in $U_{\varepsilon(g_n)} (g_n)$ for an arbitrarily large $n$ and thus $g$ is an infinite configuration.

Now additionally suppose that  $\{g_n\}$  satisfies $$\min \{d(z_1,z_2): z_1,z_2 \in   B(g_n)\} \ge s + 1/n.$$
We also  require that $\varepsilon(g_n)$ satisfies :
for any $h \in U_{\varepsilon(g_n)} (g_n)$ we have
$$\min \{d(z_1,z_2): z_1,z_2 \in B(g_n)\}  \ge s + 1/{2n}.$$
Point (2) follows directly.
\end{proof}

Fix $g \in \textit{Conf}$.  The \emph{wind-tree table} $\B^g$ is the plane $\mathbb{R}^2$ with the interiors of the union of the trees
removed.
Fix  $\theta \in \S1$.
The \emph{billiard flow $\phi_t^{g,\theta}$ in the direction $\theta$} or simply $\phi_t^\theta$  is  the free motion in the interior of $\B^g$ with elastic collision from the boundary of $\B^g$ (the boundary of the union of the trees).
Once launched in the direction $\theta$, the billiard direction can only achieve four directions $[\theta] := \{\pm \theta, \pm(\pi - \theta)\}$;
thus the \emph{phase space} $X^{g,\theta}$  
of the billiard flow in the direction $\theta$ is a subset of the cartesian product of $\B^g$ with these four directions. We agree that if a billiard orbit hits a corner of a tree, the outcome of the collision is not defined, and the billiard orbit stops there, its future is not defined anymore.
Note that in this notation  $\phi_t^\theta$, $\phi_t^{-\theta}$, $\phi_t^{\pi-\theta}$ and $\phi_t^{\theta - \pi}$
are all the same.

A flow $\psi_t$ preserving a Borel  measure $m$ is called {\em ergodic} if  for each Borel measurable set $A,  m(\psi_t(A) \triangle A) = 0 \ \forall t \in \R$ implies that $m(A) = 0$ or $m(A^c)=0$. The flow $\psi_t$ is said to have  \emph{infinite ergodic index} if for each integer 
$K \ge 1$
the  $K$-fold product flow $\psi_t \times \cdots \times \psi_t$ is ergodic with respect to the $K$-fold product measure $m \times \cdots \times m$. It is a well  known fact that in the  finite measure case the notion of   infinite ergodic index is equivalent to weak-mixing.  However we are working 
 in the context of an infinite measure preserving flow.

For each direction $\theta$, the billiard flow  $\z^{\theta}_t$ preserves 
 the area measure $\mu$ on $\B^g$ times a discrete measure on $[\theta]$, we
will also call this measure $\mu$. 
Note that $\mu$ is an infinite measure.
The billiard flow on the full phase space preserves the volume measure $\mu \times \lambda$ with $\lambda$ the length measure on $\S1$. 
Let $K\ge 1$, and let $\vt=(\theta_1,\dots,\theta_K)$ be a vector of directions. Then we note the \emph{product billiard flow} $\zt_t:=\z_t^{\theta_1}\times\cdots\times \z_t^{\theta_K}$. This flow preserves the measure $\mu^K:=\mu\times\cdots\times\mu$.
 
Now we can state our main result.
 
\begin{theorem}\label{main} For any $s>0$ there is a dense $G_\delta$ subset $G$ of $\textit{Conf}$ and a dense $G_\delta$ set of full measure of directions $\mathcal H$,  for every integer $K \ge 1$ there is a dense $G_\delta$ set  $\mathcal{H}(K)$ of full measure of $K$-tuples of directions  such that for each $g \in G$ 
\begin{enumerate}
	\item the flow  $\phi_t^{\theta}$ has infinite ergodic index for every  $\theta \in \mathcal H$ and 
	\item the flow  $\zt_t$ is ergodic for every $\vt\in \mathcal H(K)$.
\end{enumerate}
\end{theorem}
\begin{remark} We do not know that the set $\mathcal{H}(K)$ has product structure, thus (1) does not follow from (2).
\end{remark}
 
 \subsection{The precise question posed by the Ehrenfests}
Consider a large but finite number $K$ of initial points  in the wind-tree model in a given direction $\theta$.  
The Ehrenfests asked do the particles directions asymptotically equalize under the wind-tree dynamics, i.e., are there approximately
  $K/4$ particles in each direction after a large time. 
  
  This question is the motivation for our study.  
 Let
 $(\vec{z},\vec{\theta})$ denote the initial positions and velocities of these particles, and
 $f_i((\vec{z},\vec{\theta}))$ denote the number of particles pointing in the direction $i \in \{\pm \theta,\pm (\pi - \theta)\}$. 
If the functions $f_i$ were integrable, then we could give a nice answer to this question using Theorem \ref{main}, but unfortunately this is not 
the case.  We give three partial answers.
First a finite measure version.  Let $A \subset \B^g$ be a positive but finite measure subset of the wind-tree table, and let $f_i^A$ denote the function
$f_i$ restricted to the set $A \times \cdots  \times A$. This function is integrable, thus applying the Hopf ergodic theorem to the wind-tree flow yields the following corollary
(here  $K$ and $s$ are fixed, and $G$ and $\mathcal{H}$ are the dense $G_{\delta}$ sets from Theorem \ref{main})
\begin{corollary}
 For each $g \in G$, for each $A \subset \B^g$ of positive measure, for each $\theta \in \mathcal{H}$, for each $i,j$ the following limit holds almost surely as $T \to \infty$:
 $$\frac{\int_0^T f_i^A \big (\phi_t^\theta \times \cdots \phi_t^\theta (\vec{z},\vec{\theta}) \big )\, dt }{\int_0^T f_j^A \big (\phi_t^\theta \times \cdots \phi_t^\theta (\vec{z},\vec{\theta}) \big ) \, dt } \to 1.$$
 \end{corollary}
 This means that if we only count  when all the particles are in the set $A$ then the average over times of the number going in each direction is asymptotically the same.
 
If we replace the flow $\phi_t^\theta \times \cdots \phi_t^\theta$ by its first return flow $\psi_t^{A,\vec{\theta}}$ to the region $A \times \cdots \times  A$, then we can apply the Birkhoff ergodic
theorem.
\begin{corollary}
 For each $g \in G$, for each $A \subset \B^g$ of positive measure, for each $i$, the following limit holds almost surely as $T \to \infty$:
 $$\frac{1}{T} \int_0^T f_i^A \big (\psi_t^{A,\vec{\theta}} (\vec{z},\vec{\theta}) \big ) \, dt  \to \int_A f_i^A \, d\mu \times \cdots \times d\mu  = \frac{K}{4} \cdot \text{area}(A).$$
 \end{corollary}
 This means that the average over time  of the direction converges  to $K/4$, but for the first return flow.
 
 Finally we can replace the $f_i$ by integrable functions which somehow measure a similar phenomenon.  For example 
 the sum of the cubes of the reciprocal of the distance of the particles from the origin:
 $\hat{f_i}((\vec{z},\vec{\theta})) = \sum_{\{k: \vec{\theta}_k = i\}} \frac{1}{\min{(1,|z_k|^3)}}$. These functions are positive and integrable, thus we can apply the Hopf
 ergodic theorem to conclude:
 \begin{corollary}
 For each $g \in G$,  for each $\theta \in \mathcal{H}$, for each $i,j$ the following limit holds almost surely as $T \to \infty$:
 $$\frac{\int_0^T \hat{f_i} \big (\phi_t^\theta \times \cdots \phi_t^\theta (\vec{z},\vec{\theta}) \big )\, dt }{\int_0^T \hat{f_j} \big (\phi_t^\theta \times \cdots \phi_t^\theta (\vec{z},\vec{\theta}) \big )\, dt } \to 1.$$
 \end{corollary}
 This means that  then the average over time of the weighed number going in each direction is asymptotically the same.
 
 For all three results we can replace a single $\theta \in \mathcal{H}$ by a vector $\vec{\theta} \in \mathcal{H}(G)$ and state a similar result for the functions $f_i,f_i^A,\hat{f_i}$ which counts the number of particles with direction in the $i$th quadrant ($i \in \{1,2,3,4\}$). We interpret these results in the following way, if $K$ particles are launched in arbitrary generic directions, then the 
average over time of the number of particles in the different quadrants are asymptotically the same in the three senses 
mentioned above.
 
\section{Proof of wind-tree results}

\begin{figure}[h]
\begin{minipage}[ht]{0.5\linewidth}
\centering
\begin{tikzpicture}[rotate=45,scale=0.75]

\draw[thin] (0.1,0.2) rectangle +(0.5,0.5);
\draw[thin] (0.3,-0.7) rectangle +(0.5,0.5);

\draw[thin] (-0.8,-0.6) rectangle +(0.5,0.5);
\draw[thin] (-1.5,0.4) rectangle +(0.5,0.5);
\draw[thin] (1,0.4) rectangle +(0.5,0.5);

\draw[thin,lightgray](-2,-2)--(-2,2)--(2,2)--(2,-2)--(-2,-2);
 
\foreach \i in {-2.25,-1.75,...,1.75}
\draw[](\i,1.75) rectangle +(0.5,0.5);
 
 \foreach \i in {-2.25,-1.75,...,1.75}
 \draw[](\i,-2.25) rectangle +(0.5,0.5);
  
 \foreach \i in {-1.75,-1.25,...,1.25}
 \draw[](1.75,\i) rectangle +(0.5,0.5);
 
 \foreach \i in {-1.75,-1.25,...,1.25}
 \draw[](-2.25,\i) rectangle +(0.5,0.5);
\end{tikzpicture}
\end{minipage}\nolinebreak
\begin{minipage}[ht]{0.5\linewidth}
\begin{tikzpicture}[rotate=45,scale=0.75]

 \draw[thin] (0.1,0.2) rectangle +(0.5,0.5);
 \draw[thin] (0.3,-0.7) rectangle +(0.5,0.5);

 \draw[thin] (-0.8,-0.6) rectangle +(0.5,0.5);
\draw[thin] (-1.5,0.4) rectangle +(0.5,0.5);
\draw[thin] (1.1,0.4) rectangle +(0.5,0.5);
 
 \foreach \i in {-2.4,-1.85,...,2.3}
 \draw[](\i,1.95) rectangle +(0.5,0.5);
 
 \foreach \i in {-2.4,-1.85,...,2.3}
  \draw[](\i,-2.45) rectangle +(0.5,0.5);
  
 \foreach \i in {-1.9,-1.35,...,1.55}
 \draw[](2,\i) rectangle +(0.5,0.5);
 
    \foreach \i in {-1.9,-1.35,...,1.55}
 \draw[](-2.4,\i) rectangle +(0.5,0.5);

\end{tikzpicture}
\end{minipage}\nolinebreak
\caption{An 8-ringed configuration and a configuration close to it.}\label{fig}
\end{figure}
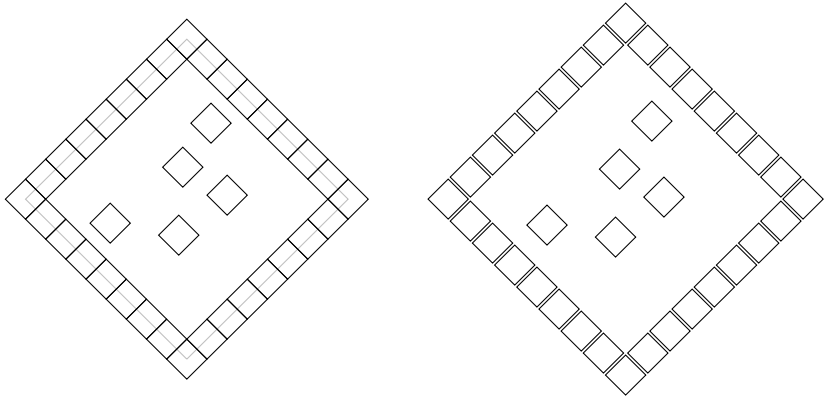
 
Fix a wind-tree configuration $g\in \textit{Conf}$. 
Fix $K$ and $n \in \N$, and  let $\B^g_{n} := \big ( \B^g \cap \{(x,y): |x| + |y| \le ns\} \big )^{K}$. 
Note that $K$ does not appear in this notation, as well as certain other notations in this section, since it is fixed throughout much of the proof.
For each $\theta$ let $[\theta]$ be the set of all possible directions under the billiard flow starting in direction $\theta$, i.e., $[\theta]=\{\pm\theta,\pm(\pi-\theta)\}$.
 Let  $[\vt] := \{ \vec \psi :  \psi_i \in [\theta_i] \text{ for all } i\}$ and 
$$X^{g,\vt}_n=\large \{(\vec z,\vec \psi): \vec z \in \B^g_n,\vec \psi \in[\vt] \large \}.$$
For each $n \ge 1$ we consider the first return flow of the product billiard flow, 
$$\z^{g,\vt,n}_t     :  X^{g,\vt}_n\to X^{g,\vt}_n.$$
For each $\theta\in\S1$, the flow $\z^{g,\vt,n}_t$ preserves the measure $\mu^K$. For sake of simplicity, we will denote $\mu^K$ by $\mu$.

\begin{proofof}{Theorem \ref{main}} 
We prove both statements with the same strategy: we choose a dense set $\{f_i\}$ of configurations which satisfy the goal dynamical property of $K$-fold ergodicity on certain compact sets.
Then we will show that wind-tree tables which are sufficiently well approximated by this dense set will satisfy 
 the  dynamical property on the whole phase space.  The proof for $K=1$ is  simpler, and we will mention the simplification in the proof even though this is
not formally necessary for the proof.

A configuration $h$ is called \emph{$n$-ringed}  if the boundary of the rhombus $\{(x,y) \in \R^2: |x| + |y| \le ns\}$ is completely covered
by trees as in Figure \ref{fig} left (i.e., the obstacles covering the boundary intersect with each other on a whole side or do not intersect at all). 

For the proof of simple ergodicity, let $\{f_i\}$ be a dense set of parameters such that each $f_i$ is an $n_i$-ringed configuration and $n_i$ is increasing with $i$. Then by \cite{KeMaSm} the billiard flow is ergodic in almost every direction inside the ring.  
So, the return flow $\z^{\vt,n}_t$  is ergodic for all $n$ such that $1\le n\le n_i$ and for almost every direction $\theta$  where $\vt$ is the vector $(\theta)$. 

Consider now the $K$-fold case.  
Let $g$ be any configuration and $\varepsilon >0 $. Let $n>\varepsilon+\frac{1}{\varepsilon }+s $. Consider the $n$-ringed configuration $f$ which  coincides with $g$   inside the the ball of radius $1/\e$ and has no additional trees. We apply Theorem 1 of \cite{MSTr3} to the table in the interior of the ring of the ringed configuration $f$ yielding a dense $G_\delta$ set $\Theta$ of full measure of directions, and  a configuration that is $n$-ringed, and is $\varepsilon$-close to $f$ and $g$
such that the flow is weakly mixing  inside the ring for all $\theta \in \Theta$.

Thus, we can find a dense set of configurations $\{f_i\}$ such that each $f_i$ is $n_i$-ringed and the flow 
 is weakly mixing for all $\theta \in \Theta$ inside the ring. Again we suppose $n_i$ is increasing with $i$.
 
Fix $K\ge 1$. Suppose that $\delta_i$ are strictly positive numbers. 
 Then the set
$$ 
\G_K := \bigcap_{m = 1}^{\infty}   \bigcup_{i=m}^{\infty} \U_{\delta_i}(f_i)
$$
is a dense $G_{\delta}$ set.  We will show that the $\delta_i$ can be chosen in such a way
that all the configurations in $\G$ are $K$-fold ergodic for all $\vt\in\Theta^K$. This $\Theta$ will be a $G_\delta$ set of full measure that has to be found in the proof. Taking intersection over $K$ will finish the proof, thus for sake of simplicity, we will fix $K$ from here on and drop it from the notations when convenient. 

Let $\{h_j\}_{j \ge 1}$ be a countable dense 
collection of continuous functions in $L^1(\R^{2K},Leb)$. For any $\vt$ and $g \in \textit{Conf}$ we think of this as a collection in $\mathcal L^1(X^{g,\vt}_n,\mu^K_n)$ in the same way as in the proof of Theorem 1 of \cite{MSTr3}.

Consider the Cesaro average
$$S^g_{n,\ell}h_j(\vec z,\vt)  := \frac{1}{\ell} \int_{0}^{\ell} h_j \big (\z^{g,\vt,n}_t(\vec z,\vt)  \big) \, dt.$$
By the Birkhoff ergodic theorem, the flow $\z^{g,\vt,n}$ is ergodic for all $n$ and for almost every $\theta$ if and only if for all $n$ and for almost all $\theta$ we have $S^g_{n,\ell}h_j(\vec z,\vt) \to \int_{X_n^{\theta}} (h_j^{\theta}(y)) \, d\mu(y)$ as $\ell$ goes to infinity for all $j \ge 1$.

Now fix $i$.  The billiard flow $\z_t^{f_i,\theta}$ is weakly-mixing inside the ring for each $\theta \in \Theta$, thus $\z_t^{f_i,\vt}$ inside the ring is ergodic for every $\vt$ in $\Theta^K$.   Thus
the first return flows  $\z_t^{f_i,\vt,n}$ are ergodic for every $\vt$ in $\Theta^K$, for all $1\le n \le n_i$.
Thus we can find  positive integers ${\ell_i} \ge n_i$,  open sets $H_i \subset \S1$ and sets $C_{n}^{f_i,\vt} \subset X_{n}^{f_i,\vt}$ so that 
$\mu(C_{n}^{f_i,\vt}) > \mu(X^{f_i,\vt}_n) - \frac1i$,
$\lambda \left (H_i\right ) > 1 - \frac1i$ and

\begin{equation}\label{eq1}
\Big |S^{f_i}_{n,\ell_i} h_j(\vec z,\vt) - \int_{X_{n}^{f_i,\vt}} h_j(y) \, d\mu(y) \Big | < \frac1i\end{equation}
for all $\vec z \in C_{n}^{f_i,\vt}$, $\vt \in (H_i)^K$,  $1 \le j\le i$, and $1 \le n \le n_i$.

Now we would like to extend these estimates to the neighborhood $\U_{\delta_i}(f_i)$
for a sufficiently small strictly positive $\delta_i$ (see Figure \ref{fig} right). For any $n$ such that $1\le n\le n_i$ let $\bar B^i_{n}$ be the intersection of $\mathcal B_{n}^{g}$ for all $g$ in the $\delta_i$-neighbourhood $\U_{\delta_i}(f_i)$. Let $\bar X_{n}^{i,\vt}:=\bar B^i_{n}\times [\vt]$.  
For every $1\le n\le n_i$ we define $\vec\psi = \vec\psi(g,f_i)$ a piecewise continuous map from $\mathcal B_{n}^{g}$ to $\mathcal B_{n}^{f_i}$. 
When convenient we will write $\vec \psi(\vec z,\vt)$ instead of  $(\vec \psi(\vec z),\vt)$.
The behavior of $\vec\psi$ will be defined coordinate by coordinate, more precisely $\vec\psi(\vec z)=(\psi(z_1)\dots,\psi(z_K))$ where $\psi$ will be defined right now. For $z$ outside the obstacles of $f_i$, we define $\psi(z)=z$. For each obstacle $O_1$ of $g$ inside the ring, we consider the corresponding obstacle $O_2$ of $f_i$ and $C_{12} = O_2\setminus O_1$. We define a direction $\xi$ that points from a corner of $O_1$ to  a corner of $O_2$ in such a way that the segment along this direction between the two corners is completely included in $C_{12}$ (as in Figure \ref{Fig2}). Then for any $z\in C_{12}$,  the image $\psi(z)$   of $z$ is the closest point in the direction $\xi$ in the table $\B_{n}^{f_i}$ (Figure \ref{Fig2}).
The difference between the Lebesgue measure of $\bar B^i_n$ and the measure of $\B_{n}^{g}$  can be made arbitrarily small by an adequate choice of $\delta_i$, simultaneously for all $g$ in $\U_{\delta_i}(f_i)$. From now on, we make the choice of $\delta_i$ such that $\mu(\bar X^{i,\vt}_{n}) = \mu(\bar B^i_n)> \mu(\B^{g}_{n}) - \frac1i$ for all $g \in \U_{\delta_i}(f_i)$ and all  $i$.

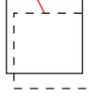
\begin{figure}[h]
\centering
\begin{tikzpicture}

\draw[thin] (0,0) rectangle +(1,1);
\draw[dashed] (0.1,-0.2) rectangle +(1,1);
\draw[red] (0.5,0.8) -- (0.4,1);

\end{tikzpicture}

\caption{The dashed obstacle is $O_1$ and the solid obstacle is the associated obstacle $O_2$. The map $\psi$ maps all the points on the red segment in the direction $\xi$ to a point in its top endpoint.}\label{Fig2}
\end{figure}

By the triangular inequality we have:
$$\begin{array}{rr}
\Big |S^g_{n,\ell_i} h_j(\vec z,\vt) - \int_{X_{n}^{g,\vt}} h_j(y) \, d\mu(y) \Big | \le& \Big |S^g_{n,\ell_i} h_j(\vec z,\vt) - S^{f_i}_{n,\ell_i}h_j(\vec \psi(\vec z),\vt)\Big | +\\
& \hspace{-2.1cm} \Big |S^{f_i}_{n,\ell_i} h_j(\vec \psi(\vec z),\vt) - \int_{X_{n}^{f_i,\vt}} (h_j(y)) d\mu(y) \Big | + \\
& \hspace{-3cm} \Big | \int_{X_{n}^{f_i,\vt}} h_j(y) \,  d\mu(y) - \int_{X_{n}^{g,\vt}} h_j(y) \, d\mu(y)\Big |
.\end{array}$$

Futhermore we choose $\delta_i$ so small that 
$$\Big | \int_{X_{n}^{f_i,\vt}\setminus \bar X_{n}^{i,\vt}} h_j(y) \,  d\mu(y) \Big |<\frac1i$$
and $$\Big | \int_{X_{n}^{g,\vt}\setminus \bar X_{n}^{i,\vt}} h_j(y) \,  d\mu(y) \Big |<\frac1i$$
thus by the triangular inequality
\begin{equation}\label{eq2}
\Big | \int_{X_{n}^{f_i,\vt}} h_j(y) \,  d\mu(y) - \int_{X_{n}^{g,\vt}} h_j(y) \, d\mu(y)\Big |<\frac2i.
\end{equation}

Now the proof bifurcates a bit according to the different cases stated in the theorem.  Consider part (1) of the theorem. So in particular $\vt$ will indicate $(\theta,\theta,\dots,\theta)$ in this part.  
Note that $\psi$ is not continuous, not invertible, and not onto. However it is not far from being continuous: $||\vec z-\vec \psi(\vec z)||_{\mathcal L^\infty}<\delta_i$ for any $\vec z\in \B^{g}_n$. 
By our convention the billiard flow stops at corners, thus any point $(\vec z,\vt)$ for which the flow is defined up to time $\ell_i$ is a point of continuity for $\z^{f_i,\vt,n}_{\ell_i}$.  Consider such a point, then  the point $\vec \psi(\z^{g,\vt,n}_{\ell_i}(\vec z,\vt)) $ stays $\delta_i$-close to $\z^{f_i,\vt,n}_{\ell_i}(\vec z,\vt)$ for $g$ in a small enough neighborhood of $f_i$; thus we can find $\delta_i>0$, an open set $\hat H_i \subset H_i$ and a set $\hat C_{n}^{i,\vt} \subset \bar X_{n}^{i,\vt}\cap  C_{n}^{f_i,\vt}$ so that that if $g \in \mathcal{U}(f_i,\delta_i)$, then
\begin{equation}\label{eq3}
\Big |S^g_{n,\ell_i} h_j(\vec z,\vt) - S^{f_i}_{n,\ell_i}h_j(\vec \psi(\vec z),\vt)\Big | < \frac2i
\end{equation}
for all $z \in \hat C_{n}^{i,\vt}$, $\theta \in \hat H_i$ (here $\vt = (\theta, \theta, \dots , \theta)$),  $1 \le j\le i$, $1 \le n \le n_i$; and $\mu(\hat C_{n}^{i,\vt}) > \mu(\B^g_n)-\textstyle\frac2i$ and $\hat H_i$ is of measure larger than $1-\textstyle\frac2i$.

Since $\lambda(\hat H_i) > 1 - 2/i$, the set $
\mathcal{H} = \cap_{M=1}^\infty \cup_{i=M}^\infty \hat H_i$ has full measure. 
Fix $g \in \G$ and  $\theta \in \mathcal{H}$, then there is an infinite sequence $i_k$ such that $g \in \U_{\delta_{i_k}}(f_{i_k})$ and $\theta \in \hat H_{i_k}$.  
Fix $n \ge 1$ and consider
$\mathcal{C}^{g,\vt}_{n}  := \cap_{M=1}^\infty \cup_{k=M}^\infty \widehat C^{i_k,\vt}_{n}$.
Recall that we made the choice of $\delta_i$ such that $\mu(\bar X^{i_k,\vt}_{n}) > \mu(\B^{g}_{n}) - \frac1{i_k}$. Since $\mu(\widehat C^{i_k,\vt}_{n}) > \mu(\bar X^{i_k,\vt}_{n}) - \frac1{i_k}$, it follows that  $\mu(\mathcal{C}^{g,\vt}_{n}) = \mu(\B^{g}_{n})$. 

Suppose $g \in \G$. 
Thus  for $\theta \in \mathcal{H}$, for each $n \ge 1$ 
the three inequalities \eqref{eq1}, \eqref{eq2}, \eqref{eq3} imply that
$$
|S^g_{n,\ell_{i_k}}(h^\theta_j) - \int_{X_n^{\theta}} (h_j^{\theta}(\vec z,\vt)) d\mu|<\frac5i
$$
for all $z \in \hat C_{n}^{i,\vt}$, $\theta \in \hat H_i$ (here $\vt = (\theta, \theta, \dots , \theta)$),  $1 \le j\le i$, $1 \le n \le n_i$
and thus 
\begin{equation}\label{e11}
\lim_{k \to \infty} S^g_{n,\ell_{i_k}}(h^\theta_j) \to 
\int_{X_n^{\theta}} (h_j^{\theta}(\vec z,\vt)) d\mu
\end{equation}
for all $(\vec z,\vt)$ in $\mathcal{C}^{g,\vt}_{n}$, for each $j \ge 1$.
The $h^{\theta}_j$ are dense in  $L^1(X_n^\theta,\mu)$ and $\lim_{k \to \infty} \ell_{i_k} = \infty$, thus Equation \eqref{e11} together
with the Birkhoff ergodic theorem imply that for each $n \ge 1$, the first return flow $\z^{g,\vt,n}_{t}$ is 
ergodic for all $\theta \in \mathcal{H}$. This implies the ergodicity of the billiard flow $\phi_t^{g,\vt}$ in every direction in $\mathcal{H}$.

For part (2) of the theorem we have to slightly modify the previous arguments.  The only difference being that  the set of directions we construct depends
on $K$.
For any point $(\vec z,\vt)$ of continuity of $\z^{f_i,\vt,n}_{\ell_i}$, the point $\z^{g,\vt,n}_{\ell_i}(\vec z,\vt) $ varies continuously with $g$ in a small neighborhood of $f_i$; thus 
we can find $\delta_i>0$, an open set $\hat H_i(K) \subset (H_i)^K$ and a set 
$\hat C_{n}^{i,\vt} \subset \bar X_{n}^{i,\vt}\cap  C_{n}^{f_i,\vt}$
 so that that if $g \in \mathcal{U}(f_i,\delta_i)$, then
$$\Big |S^g_{n,\ell_i} h_j(\vec z,\vt) - S^{f_i}_{n,\ell_i} h_j(\vec \psi(\vec z),\vt)| < \frac2i$$
for all $z \in \hat C_{n}^{i,\vt}$, $\vt \in \hat H_i(K)$, $1 \le n \le n_i$, $1 \le j\le i$; and $\mu(\hat C_{n}^{i,\vt}) >  \mu(\B^g_n)-\textstyle\frac2i$ and $\hat H_i(K)$ is of measure larger than $1-\textstyle\frac2i$.

Since $\lambda(\hat H_i) > 1 - 2/i$, the $G_\delta$ set $
\mathcal{H}(K) = \cap_{M=1}^\infty \cup_{i=M}^\infty \hat H_i(K)$ has full measure.
The rest of the proof of part (2) is identical to that of part (1).
\end{proofof}

\subsection{Generalization}
If we consider a subset $C$ of $(\textit{Conf},d_H)$ which is itself a Baire set  such that
the set $\{h: \nolinebreak h \text{ is } $N$-\text{ringed for } N \ge N_0\}$  is dense in $C$ for each $N_0 \ge 1$ then Theorem  \ref{main} holds in $(C,d_H)$ as well.
In particular the set of configurations considered in the articles \cite{MSTr1},\cite{MSTr2} is a Baire subset of $(\textit{Conf},d_H)$
thus Theorem \ref{main} holds in that context as well.

\section{Appendix}

\begin{proofof}{Proposition  \ref{p1}}
Let $(g_i)_{i \in \N}$ be a sequence of configurations. Consider $\varepsilon_n = 1/n$. 
Let $B_\varepsilon := \{ z: \rho(z,\infty) > \varepsilon) \}$.
Let $k_j$ be the cardinality of $g_j \cap B_{\varepsilon_1}$.
The sequence $k_i$ only take a finite number of values. Thus we can choose subsequence $(g_{j})_{j \in J_0}$ such that  the sequence $(k_{j}: j \in J_0)$ is constant, call this constant $c_1$. 

If $c_1 = 0$ then for each $j \in J_0$ let $g^1_j$ be the empty configuration. Otherwise for each $j \in J_0$  let
 $g^{1}_j := \{z^1_{1,j},\dots, z^1_{c_1,j}\} = g_{j} \cap B_{\varepsilon_1} $.  For each $j$ we think of $g^1_j$ as
a finite configuration in $\textit{Conf}$, but also as a point in $\mathbb{S}^{c_1}$.  By compactness of 
$ \mathbb S^{c_1}$ 
we can find a subsequence $J_1 \subset J_0$ such that the $(g^1_j: j \in J_1)$ converge to a point $g^1 := (z^1_1,\dots,z^1_{c_1})\in \mathbb{S}^{c_1}$. Note that $d(z^1_i,z^1_j) \ge s$ for all $i \ne j$, thus $g^1 \in \textit{Conf}$.  Furthermore we have
$g^1_j \in \U_{\varepsilon_1}(g^1)$ for all sufficiently large $j \in J_1$. 
Repeat this argument for $n=2$ to produce  a  subsequence $J_2 \subset J_1$
which converge 
to a point $g^{2} \in \mathbb{S}^{c_2}$. Again we have $g^2 \in \textit{Conf}$ and $g^2_j \in \U_{\varepsilon_2}(g^2)$ for all sufficiently large $j \in J_2$. Note that $c_2 \ge c_1 \ge 0$ and for $j=1,\dots,c_1$ we have $z^2_j = z^1_j$.
Repeat this construction for each $n$.  Finally we define $g$ to be an almost countable collection of points such that
every $z \in g$ is in $g^k$ for all sufficiently large $k$. By construction for any $z_1,z_2 \in g$ we have $d(z_1,z_2) \ge s$,
thus $g \in \textit{Conf}$.
Note that  $g$ can be  an infinite, finite, or even the empty configuration.

We claim that $g$ is an accumulation point of the sequence $(g_j)_{j \in \N}$.  
Fix a neighborhood $\U$ of $g$.  Choose $n$ so large that ${\U_{\varepsilon_n}(g)} \subset \U$.  
By construction of $g$  and $g^n$ we have $\U_{\varepsilon_n}(g^n) = \U_{\varepsilon_n}(g)$.
The result follows since $g_j \in \U_{\varepsilon_n}(g^n)$ for all sufficiently large $j \in J_n \subset \N$.
\end{proofof}

\begin{remark}
If we remove the empty and finite configurations from $\textit{Conf}$ the space is not even locally compact.
\end{remark}

\end{document}